# New approach for Finite Difference Method for Thermal Analysis of Passive Solar Systems

Stanko Shtrakov and Anton Stoilov

*Department of Computer's systems, South - West University "Neofit Rilski", Blagoevgrad, BULGARIA,*
*(Dated: February 17, 2005)*

Mathematical treatment of massive wall systems is a useful tool for investigation of these solar applications. The objectives of this work are to develop (and validate) a numerical solution model for predication the thermal behaviour of passive solar systems with massive wall, to improve knowledge of using indirect passive solar systems and assess its energy efficiency according to climatic conditions in Bulgaria. The problem of passive solar systems with massive walls is modelled by thermal and mass transfer equations. As a boundary conditions for the mathematical problem are used equations, which describe influence of weather data and constructive parameters of building on the thermal performance of the passive system. The mathematical model is solved by means of finite-differences method and improved solution procedure. In article are presented results of theoretical and experimental study for developing and validating a numerical solution model for predication the thermal behaviour of passive solar systems with massive wall.

## I. INTRODUCTION

The concept of passive solar systems is well-known method for use of solar energy as a source of heating in buildings. There is a vast literature on this technology, but real objects (houses) with passive solar systems are still rare. The major impediments to increase market penetration of passive solar systems is the lack of available information and experience data for the efficiency and constructive parameters of passive solar elements.

The main concept of indirect passive solar systems is Trombe-Michel wall. Most experimental and theoretical data, published on Trombe wall performance are in form of overall building performance. Data of overall performance is of limited use, as it only provides seasonal estimates of heat gains for specific building designs, wall patterns and climates. Because of the large number of parameters and the wide range of weather conditions, which influence the operation of massive walls, the assessment of the thermal behavior requires the use of thermal simulation techniques.

Literature review shows that the problem of passive solar systems with massive walls is ordinarily modelled by thermal and mass transfer equations [4,7]. As boundary conditions for the mathematical problem must be used equations, which describe influence of weather data and constructive parameters of building on the thermal performance of the passive system. The mathematical model, composed for the massive wall performance, is usually very complicated and for solving the mathematical system of equations it is necessary to apply a different set of assumptions.

The purpose of this article is to present the results of theoretical and experimental study for developing and validating a numerical solution model for predication the thermal behaviour of passive solar systems with massive wall.

## 2. MATHEMATICAL MODEL

The simulation scheme of typical passive solar system with massive wall is shown in Fig. 1. The massive wall is usually mounted on the south facade of the house. It comprises three layers: a transparent cover (one or two glasses or plastic plates), a massive wall (masonry, concrete) and air gap between transparent cover and massive wall. At the bottom and the top of the massive wall, there are vents for allowing an air circulation between air gap and room space. The external transparent cover transmits solar radiation in, but holds back heat. Wall surface is painted black at its outer side and act as an absorber of solar radiation. It stores heat from the day and releases it with time to the room space by radiative heat transfer. The air layer between transparent cover and massive wall is heated in contact with the wall surface, rises and circulates towards the room (when the vents are open).

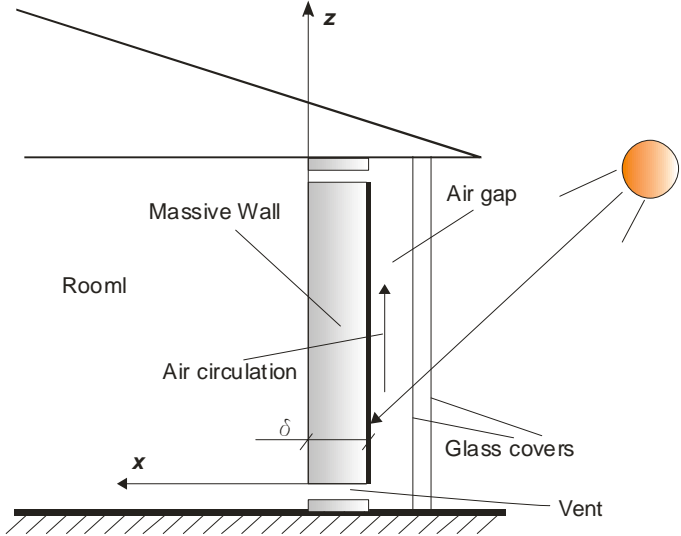

Fig.1. Scheme of 'Trombe Wall' System

In summer, inhabitants close the vents and air circulation is stopped during the day. During the night, vents in upper end of glazing can be opened and natural ventilation in room space can be organized.

The thermal analysis of such system is very complicated. Mathematical model is based on the transient performance of the system. It comprises energy balance equations, written for each element of the system. Since the wall is taken to be large (in comparison with wall thickness), the temperature variation in **y** – direction (wall width) will be neglected, and only two-dimensional problem can be considered - the hight (z – direction) and thickness (x – direction). The governing energy conservation equation of heat transfer in massive walls is:

$$\frac{\partial T}{\partial \tau} = a \left( \frac{\partial^2 T}{\partial x^2} + \frac{\partial^2 T}{\partial z^2} \right) \tag{1}$$

where **T** refers to the temperature in the wall, $\tau$ is time variable, $x$ and $z$ – space variables and $a$ is a material constant. The boundary conditions, needed for the solution of eqn (1), are derived from an energy balance for elements of passive system. The cover glazing is assumed to absorb no solar radiation and only exchanges heat by convection and radiation with the wall surface and the ambient. Heat transfer rate in elements of the massive wall is given by:

at $x = \delta$ (an inner wall surface)

$$\lambda \frac{dT}{dx} = h_{cRoom}(T_{w0} - T_r) + h_{rRoom}(T_{w0} - T_{wall}) \tag{2}$$

at $x = 0$: (an outer wall surface)

$$-\lambda \frac{dT}{dx} = q_s(\tau\alpha)_e + h_{cgap}(T_{ag} - T_{wn}) + h_{rgap}(T_{g2} - T_{wn}) \tag{3}$$

at inner glass cover:

$$h_{cgap}(T_{ag} - T_{g2}) + h_{rgap}(T_{w0} - T_{g2}) = h_{12}(T_{g2} - T_{g1}) \tag{4}$$

at outer glass cover:

$$h_{c\infty}(T_{g1} - T_a) + h_{r\infty}(T_{g1} - T_{sky}) = h_{12}(T_{g2} - T_{g1}) , \tag{5}$$

where $T_{wn}, T_{w0}$ are temperatures of wall surfaces ($x=\delta$ and $x=0$)

$T_{ag}, T_{g1}, T_{g2}, T_a, T_{wall}, T_r, T_{sky}$ - temperatures of air in the gap, glass covers, ambient air, averages of room's walls, air in room and sky, respectively

$h_{c..}$ - Convective transfer coefficients, W/m$^{2o}$K,

$h_{r..}$ - Radiant transfer coefficients, W/m$^{2o}$K,

$h_{12}$ - heat transfer coefficient in space between glass covers, W/m$^{2o}$K.

Next the air gap is considered. A differential formulation, which includes terms due to thermal capacity and convective heat transfer from the wall surface and glazing cover to the air, leads to the equation:

$$\rho G_{ag} c_p \frac{dT_{ag}}{dz} = (h_{cgap}(T_{w0} - T_{ag}) + h_{cgap}(T_{g2} - T_{ag}))Bdz \tag{6}$$

where $\rho$ is the air density [kg/m$^3$]; $c_p$ - heat capacity of air [J/kg$^o$K]; $G_{ag}$ - air flow rate in a gap, [m$^3$/s]; $B$ - wall width [m], $(\tau\alpha)_e$ - an absorptance-transmittance product for the total insolation on a vertical surface

The air flow rate for air circulation in the air cavity is determined by the average air velocity. According to J.A.Duffie and W.A. Beckman [4], natural convection in cavity can be assessed by next expression:

$$V = \sqrt{\frac{2gH}{C_1(A_g/A_v)^2 + C_2} \frac{T_m - T_r}{T_m}} \qquad (7)$$

where $C_1$ and $C_2$ are constants that depend on hydraulic characteristics of the gap, $A_g, A_v$ - gap area and vent area [m$^2$] ; $H$ - wall height, $T_m, T_r$ - average temperature of air in gap and temperature of air in room.

System (1) – (7) is unsteady two-dimensional mathematical model of massive wall system. The model has a combined system of algebraic and differential equations as a boundary condition. This model is rather difficult to solve due to the complicate thermal and mass transfer processes in system. The main problems which arise during solving processes are:

- Temperature variation in direction $z$ (according to the wall height) requires solving all equations (2) – (7) with respect to this variation. Hence, it is necessary to consider vertical temperature variations in the glasses and gap. This can be done only by including new heat transfer equations (such as equation (1)) for glass covers. This temperature variation is caused mainly by air circulation in the cavity (air convection).

- The circulation in the air cavity depends on the value of air velocity and therefore, on the temperature difference in down and upper part of the air gap - eqn. (7). At the same time, the temperature difference depends on the convective transfer coefficient $h_{cgap}$, which is function of the air velocity. This determines the mathematical model as nonlinearly and requires special iteration procedures for solving the heat and mass transfer equations.

- The complicated form of equations for boundary conditions presumes many difficulties in trying to solve the mathematical equations by regular numerical scheme. This goal is complicated additionally by unregulated variation of ambient climatic parameters (ambient temperature and solar radiation).

From the point of view of engineering application, for the simulation model of passive solar system with massive wall the following assumptions can be grounded:

- The literature review [1,2,5] and calculations, we carried out, showed that the heat transfer by convection with air circulation is up to 10 - 15% of all heat transfer in massive wall. This heat transfer determines very small vertical temperature variation in the wall (0.2 - 0.5 $^o$C), because of the temperature equalization by heat conductivity and the large heat capacity of the wall. On the base of these results, the model can be simplified to one-dimensional one by assuming the different layers of wall construction to be at uniform temperature at any given time. In this way, equation (1) can be simplified to one-dimensional problem, referred only to variable $x$.

- Thermal characteristics of the massive wall and the air flow in cavity are considered as constants, because of small temperature variation in thermal and mass transfer processes. Convection and radiation heat transfer coefficients in energy balance equations are treated as depending on velocity and temperature difference between corresponding elements.

- Lastly, an iterative calculation process can be organized, if equation (6) has solved separately, by using uniform temperatures of wall surface and glass cover (a numerical method for solving differential equation can be used with considering these temperatures known). The solution will determine the air temperature variation in the air cavity. Separate finite difference method in $z$ direction has been used for solving the equation (6). On the base of received temperature rise in the duct, equation (6) can be rewritten as algebraic scheme using mean air temperature in the air duct:

$$\rho G_{ag} c_p \frac{T_m - T_r}{HB} = h_{cgap}(T_{wn} - T_{ma}) + h_{cgap}(T_{g2} - T_{ma}) \qquad (8)$$

where
$$T_{ma} = \frac{T_m + T_r}{2}$$

and $T_m$ and $T_\gamma$ are the air temperatures in bottom and top of the air gap. After solving the main task (1), (2), (3), (4), (5) and (8) and receiving values for temperatures of wall surface and inner glass, the equation (6) can be solved again and all procedure can be repeated. Process can be continued until sufficient accuracy is arrived.

## 3. FINITE DIFFERENCE APPROXIMATION

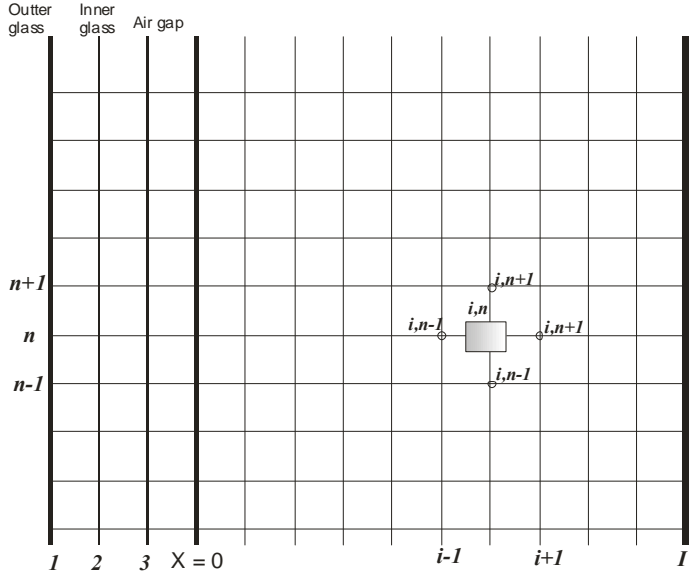

Fig.2. The mesh in time and space

The finite-difference form of differential equation (1) is derived by integration over control volume surrounded the typical node *i, n* in solution grid (Fig.2). The indexes *i* and *n* refer to the thickness (*x*) and the time ($\tau$) variable, respectively. An implicit time approximation, which is stable for forward integration in time, is developed for transient differential equations. In this case a set of simultaneous equations needs to be solved at each time step. If the time interval is named $\Delta\tau = (\tau_n, \tau_n+\Delta\tau) = (\tau_n, \tau_{n+1})$, the time derivative can be written using forward Euler formula for discretization:

$$\frac{\partial T}{\partial \tau} = \frac{T^{n+1}_i - T^n_i}{\Delta \tau} \quad (9)$$

For the space derivative is applied symmetrical Cranc-Nicolson's scheme [3] for discretization:

$$\frac{\partial^2 T}{\partial x^2} = \sigma \frac{T^{n+1}_{i+1} - 2T^{n+1}_i + T^{n+1}_{i-1}}{\Delta x^2} + (1-\sigma)\frac{T^n_{i+1} - 2T^n_i + T^n_{i-1}}{\Delta x^2} \quad (10)$$

where $\sigma$ is a weight coefficient.

After substituting (9) and (10) in (1) and rearranging, the following general approximation is received:

$$T^{n+1}_{i-1} - (2 + \frac{\nabla x^2}{\sigma \nabla \tau\, a^n_i})T^{n+1}_i + T^{n+1}_{i+1} = F^n_i, \quad i = 1,2,..I \quad (11)$$

where: $\quad F^n_i = \frac{1-\sigma}{\sigma}(T^n_{i-1} - 2T^n_i + T^n_{i+1}) + \frac{\Delta x^2}{\sigma \Delta\tau\, a^n_i} T^n_i$

This is a system of (*I*-2) algebraic equations with *I* (*i* = 1,2....*I*) unknown node temperatures (with upper index *n+1*). Temperatures with upper index *n* are considered as known, which are received from calculation, made in former time step or from initial conditions in the first time step. Equation (11) yields a special diagonal matrix of equations for the time layer *n+1* at each space point. Boundary conditions (2), (3), (4), (5) and (8) must be added to complete the system. These equations include a new unknown temperatures (*T_{ag}*, *T_{g1}*, *T_{g2}*), which requires preliminarily solving the boundary equations system.

Because of complicated nature of boundary conditions and numerous difficulties, which appear in solution process, a new procedure for completing the algebraic system of equations is proposed. The mesh is expanded by adding three new space layers, numbered 1, 2 and 3, as it is shown on Fig.2. These layers correspond to the elements of Trombe wall - two glass covers and air in the gap with temperatures *T_{ag}*, *T_{g1}*, *T_{g2}*. This means that, to the algebraic system (11), it is necessary to add algebraic equations (3), (4) and (8). As boundary conditions in this new system must be used only equations (2) and (5). In this way, mathematical task becomes considerably easier to solve, because of simplified boundary condition system.

## 4. SOLUTION PROCEDURE

Equations (11) can be solved by standard algebraic methods. Because of simple form of algebraic system (11), the well-known procedure with twofold calculation passage in the space direction of mesh is used [8]. This method is applicable for algebraic systems in form as follow:

$$A_i T_{i-1} + B_i T_i + C_i T_{i+1} = F_i \quad (12)$$

where $i = 1,2........I-1$, whit boundary conditions:

$$T_0 = a_0 T_1 + b_0 \quad \text{and} \quad T_I = a_I T_{I-1} + b_I. \tag{13}$$

This algebraic system is similar to our problem (11), (2), (3), (5) and (8) with respect to the unknown temperatures (superscript index $n+1$). In the next considerations the upper indexes of variables in finite difference equations can be omitted for simplicity.
Solution for above system is wanted in form as follow:

$$T_i = \alpha_i T_{i+1} + \beta_i, \tag{14}$$

where $\alpha_i$ and $\beta_i$ are unknown coefficients. This equation can be written for all indexes $i$, including $i-1$:

$$T_{i-1} = \alpha_{i-1} T_i + \beta_{i-1} \tag{15}$$

Substituting (14) in (12) and rearranging, it can be received next formulas for coefficients:

$$\alpha_i = \frac{C_i}{A_i \alpha_{i-1} + B_i}, \quad \beta_i = \frac{F_i - A_i \beta_{i-1}}{A_i \alpha_{i-1} + B_i} \tag{16}$$

If it is take into account that for $i = 0$: $\alpha_0 = a_0$ and $\beta_0 = b_0$ (from boundary condition eqn.(13)), coefficients $\alpha_i$ and $\beta_i$ can be calculated by using recurrent equation (16). This is the first calculation passage in the space direction of the solution grid. Two equations are available for the last node $I$ of the grid: second part of equation (13) and equation (14) for $I-1$ node: $T_{I-1} = \alpha_{I-1} T_I + \beta_{I-1}$.
From these two equations, it is possible to determine $T_I$:

$$T_I = \frac{a_I \beta_{I-1} + b_I}{1 - a_I \alpha_{I-1}} \tag{17}$$

Knowing $T_I$, it can be made second calculation passage through the solution grid by using recurrent equation (14) to calculate temperatures $T_i$ ($i = I-1, I-2......0$) of all nodes of the grid. This is a ordinarily procedure for solving the problem in time step. Receiving the temperatures for time step $n+1$, it is possible to make next time step.

This simple procedure for solving the algebraic system of finite difference approximation can not be used directly for the mathematical problem of massive wall system described above. System of equations (2), (3), (4) and (8), which is added to main system (11), is not fully compatible with the system (12). Equation (4) comprises four unknown temperatures - two temperatures of glass covers $T_{g1}$ and $T_{g2}$, temperature of air in gap $T_{ag}$ and surface temperature $T_{wn}$. This means, that equation (4) can not be expressed in form as the equation (14).

To use a technique similar to mentioned above procedure for solving equations (12), a modified procedure can be used. Instead of eqn (15), a new temperature function (with three consecutive temperatures) can be considered:

$$T_i = \alpha_i T_{i+1} + \gamma_i T_{i+2} + \beta_i \tag{18}$$

This means that, tree unknown coefficients $\alpha_i$, $\beta_i$ and $\gamma_i$, must be calculated in first calculating passage. After doing similar steps as it has been made above (equation 16), equations for unknown coefficients $\alpha_i$, $\beta_i$ and $\gamma_i$ can be easier to receive. Here, these coefficients will not be described for the common case, but for the special case of mathematical model for passive solar system with massive wall. After rearranging equations of boundary conditions in standard form (equation (12)), with respect to temperature function (18) and indexes in numerical grid (fig.2), the following equations for coefficient $\alpha$, $\beta$ and $\gamma$ have been received:
- **for outer glass cover** – equation (5). It can be rewritten in form like eqn.13:

$$T_{g1} = \alpha_1 \cdot T_{g2} + \beta_1, \tag{19}$$

where the coefficients with index $i=1$ are:

$$\alpha_1 = \frac{h_{12}}{h_{12} + h_{c\infty}}, \quad \beta_1 = \frac{h_{c\infty}}{h_{12} + h_{c\infty}} T_a, \quad \gamma_1 = 0 \tag{20}$$

- **for inner glass cover** – equation (4). Substituting $T_{g1}$ from (19) and rearranging, the equation (4) can be transformed in form as eqn.18 with coeficiens ($i=2$):

$$\alpha_2 = \frac{h_{cgap}}{h_{12}(1-\alpha_1) + h_{gap}}, \quad \beta_2 = \frac{h_{12}\beta_1}{h_{12}(1-\alpha_1) + h_{gap}} \quad \gamma_2 = \frac{h_{rgap}}{h_{12}(1-\alpha_1) + h_{gap}}$$

(21)

where $h_{gap} = h_{cgap} + h_{rgap}$

- **for air gap** – equation (8). Using analogous transformation, it can also be received equation in form as eqn.18 with coefficients (*i=3*):

$$\alpha_3 = \frac{BHh_{cgap}(\gamma_2 + 1)}{D_x}, \quad \beta_3 = \frac{2\rho G c_p BT_r + BHh_{cgap}\beta_2}{Dx}, \quad \gamma_3 = 0 \qquad (22)$$

where $D_x = BH\,h_{cgap}(\gamma_2+1) + 2\,\rho\,G\,c_p\,B$

- **for wall surface** – equation (2). Here is necessary to approximate the x direction derivative of temperature by finite difference. The appropriate coefficients for algebraic equation are: (*i=4*):

$$\alpha_4 = \frac{\lambda\sigma/\Delta x}{E_x}, \quad \beta_4 = \frac{F_x - \sigma h_{rgap}\alpha_2\beta_3 - \sigma h_{rgap}\beta_2 - \sigma h_{cgap}\beta_3}{E_x}, \quad \gamma_4 = 0 \qquad (23)$$

where:

$$E_x = \sigma\,h_{rgap}\,\alpha_2\,\alpha_3 + \sigma\,h_{rgap}\,\gamma_2 + \sigma\,h_{cgap}\,\alpha_3 + \sigma\,h_{cgap} - \frac{\lambda\sigma}{\Delta x} - \frac{\lambda\Delta x}{2\,a_4^n\,\Delta\tau}$$

$$F_x = \lambda(1-\sigma)\frac{T_5^n - T_4^n}{\Delta x} + \lambda\frac{\Delta x\,T_4^n}{2\,a_4^n\,\Delta\tau} - q_s(\tau\alpha)_e - (1-\sigma)h_{cgap}(T_3^n - T_4^n) + h_{rgap}(T_2^n - T_4^n)$$

- **for ordinary wall layer** - Here is valid the standard transformation (16) for equation (11) (*i*):

$$\alpha_i = \frac{1}{2 + \frac{\Delta x^2}{\sigma\Delta\tau\,a_{i-1}^n} - \alpha_{i-1}}, \quad \beta_i = (\beta_{i-1} - F_i^n)\alpha_i, \quad \gamma_i = 0 \qquad (24)$$

Temperature of the last wall's layer (inner wall surface) can be defined by following equation:

$$T_I^{n+1} = \frac{\frac{\sigma\lambda}{\Delta}x\,\beta_{I-1} - FF}{\sigma\,h_{ROOM} + \frac{\sigma\lambda}{\Delta x} + \frac{\rho c_p \Delta x}{\Delta\tau} - \frac{\sigma\lambda}{\Delta x}\alpha_{I-1}} \qquad (25)$$

where

$$FF = \sigma h_{ROOM}T_r + (1-\sigma)(T_{I+3}^n - T_r)h_{ROOM} + \frac{(1-\sigma)\lambda}{\Delta x}(T_{I+3}^n - T_{I+2}^n) - \rho c_p\frac{\Delta x}{\Delta\tau}T_{I+3}^n$$

and $h_{ROOM} = h_{cROOM} + h_{rROOM}$

Knowing temperature $T_I^{n+1}$, other node temperatures can be calculated with the recurrent formula (18). It is possible, because coefficient $\gamma$ for ordinary wall layers is zero. This is the second calculation passage on the grid.

To solve the mathematical problem by proposed algorithmic scheme, the heat transfer coefficients and the air velocity in gap, must be known in advance. These can be calculated with regarding temperatures and other needed variables from previous time-step calculation or from initial conditions in starting the calculations. For improving the precision of numeric calculations, an iterative calculation have been organised until needed accuracy has arrived.

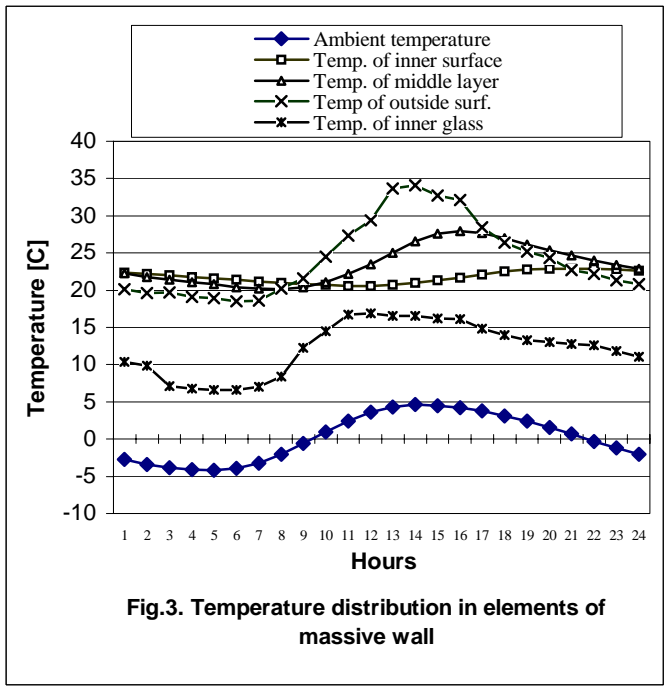

Fig.3. Temperature distribution in elements of massive wall

## 5. NUMERICAL EXAMPLES

To verify the applicability of the above-proposed technique, a large number of numerical examples have been carried. In Fig.3 and 4 is shown example of calculations for passive solar system with massive wall. Climatic data are for Sofia, Bulgaria. Passive system is with south facing concrete wall with dimensions: height - 3 m; width - 3.5 m and thickness - 0.3m. Climatic data (solar radiation and ambient temperature) are for February. The set of month's daily distribution of solar radiation and ambient temperature, estimated in hour-by-hour period have been used. Five days period of simulation calculations was needed to exclude influence of initial conditions. Following parameters are shown in Fig. 3 and 4: ambient temperature $T_a$, solar radiation $q_s$, inner glass temperature $T_{g2}$, outer surface temperature $T_{wn}$, wall temperature in middle layer $T_{wm}$, inner surface temperature $T_{w0}$, heat losses to the ambient air $q_t$, heat transferred by convection $q_f$, heat transferred by radiation $q_r$.

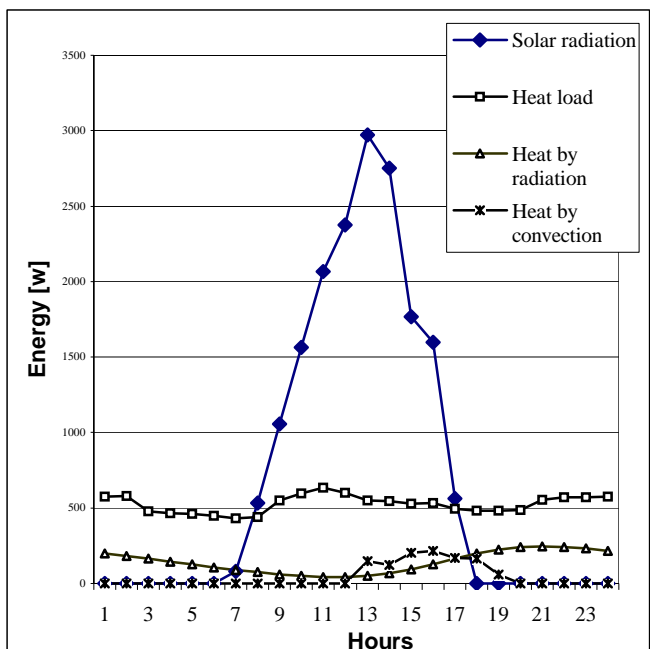

Fig.4. Energy balance in massive vall

## 6. CONCLUSIONS

In this paper, has been presented a mathematical model for passive solar system with massive wall. A finite difference solution scheme, based on implicit method was developed for solving the combined system of algebraic and differential equations. A new solution procedure, suitable for the presented mathematical model was suggested.

Computer program for simulation calculations was created, and with extensive numerical experiments, the applicability of presented model was verified.

The results from presented mathematical model would help researches in field of passive solar systems to increase their knowledge and experience for thermal and mass transfer processes. The designers of passive solar systems can use this model to select optimal constructive parameters of the massive wall.